\documentclass[12pt]{amsart}
\usepackage{amsmath,amssymb}
\usepackage[mathscr]{eucal}
\begin{document}

\newcommand{\bbR}{{\mathbb R}}
\newcommand{\bbN}{{\mathbb N}}
\newcommand{\bbI}{{\mathbb I}}
\newcommand{\bbQ}{{\mathbb Q}}
\newcommand{\bbP}{{\mathbb P}}
\def\A{\mathscr{A}}
\def\B{\mathscr{B}}
\def\C{\mathscr{C}}
\def\D{\mathscr{D}}
\def\E{\mathscr{E}}
\def\F{\mathscr{F}}
\def\H{\mathscr{H}}
\def\K{\mathscr{K}}
\def\L{\mathscr{L}}
\def\M{\mathscr{M}}
\def\N{\mathscr{N}}
\def\G{\mathscr{G}}
\def\P{\mathscr{P}}
\def\U{\mathscr{U}}
\def\V{\mathscr{V}}
\def\W{\mathscr{W}}
\renewcommand{\ss}{\subseteq}
\newcommand{\dsize}{\displaystyle}
\def\dom{\mathop{\mathrm{dom}}}
\newtheorem{thm}{Theorem}[section]
\newtheorem{cor}[thm]{Corollary}
\newtheorem{lem}[thm]{Lemma}
\newtheorem{prop}[thm]{Proposition}
\newtheorem{remark}[thm]{Remark}
\newtheorem{ex}[thm]{Example}
\font\sans=cmss10
\def\axiom#1{{\sans (#1)}}
\def\wk{\text{weak}}
\def\domc#1{{\dom c_{#1}}}
\def\sm{\smallsetminus}

\title{A remark on a theorem of Iliadis concerning isometrically containing mappings}
\author{El\.{z}bieta Pol and Roman Pol}
\address{University of Warsaw \and  University of Warsaw}
\email{E.Pol@mimuw.edu.pl \and  R.Pol@mimuw.edu.pl}
\date{\today}

\begin{abstract}
In this paper we give an alernative proof and a refinement of a recent result of S.D.Iliadis, concerning isometrically containing mappings. We address also a recent related result by A.I.Oblakova.
\end{abstract}

\maketitle


\section{Introduction}

Following Stavros Iliadis \cite{I4}, we shall say that a continuous surjection $\;F : X \to Y\;$ between separable metric spaces is isometrically containing for a class $\;\F\;$ of continuous surjections between compact metric spaces if for each $\;f \in \F\;$ there are isometric embeddings \\$i : dom f \to X$ and $j : ran f \to Y$ of the domain and the range of $f$, respectively, such that $F \circ i = j \circ f$.

Iliadis \cite{I4} proved the following theorem.

\bigskip

{\bf  Theorem 1} (S.D. Iliadis). {\it For any countable ordinals $\alpha , \beta$ there is a continuous surjection $F : X \to Y$ between complete separable metric spaces with small transfinite dimensions {\rm ind}$X = \alpha$, {\rm ind}$Y = \beta$, which is isometrically containing for the class $\F (\alpha , \beta )$ of continuous surjections $f$ between compact metric spaces with {\rm ind}$(dom f) \leq \alpha$ and {\rm ind}$(ran f) \leq \beta$.}

\bigskip

The proof given by Iliadis is based on a rather refined method, developed in his earlier works \cite{I1} - \cite{I3}.

The aim of this remark is to provide another proof of this theorem, using some results from \cite{P1}, cf. also \cite{PP} and \cite{P2}.

This approach gives also some refinement of the Iliadis theorem to the following effect (recall that the dimension of a mapping is the supremum of dimensions of its fibers):

\bigskip

{\bf Theorem 2.} {\it Given countable ordinals $\alpha$, $\beta$  and an $n$ which is a natural number or $\infty$, there is a continuous surjection $F : X \to Y$ between complete separable metric spaces with {\rm ind}$X = \alpha$, {\rm ind}$Y = \beta$, such that {\rm dim}$F \leq n$ and $F$ is isometrically containing for the class $\F (\alpha , \beta , n )$ of the maps $f \in \F (\alpha , \beta )$ with {\rm dim}$f \leq n$.}

\section{Proof of Theorem 2}

 In the sequel, we shall adopt the notation from the proof of Proposition 5.1.1 in \cite{PP}.

Let $M$ be a complete separable metric space isometrically universal for separable metric spaces, let $\K (M)$ be the hyperspace of compact subsets of $M$ equipped with the Vietoris topology, and for a countable ordinal $\gamma$, let $\K _{\gamma} (M)= \{ K \in \K (M): \;{\rm ind}K \leq \gamma \}$.

We shall consider also $M$ as a topological subspace of the Hilbert cube $I ^{\infty}$, and it should be clear from the context, when we refer to the fixed universal metric on $M$ (which does not extend over $I^{\infty}$) or we just deal with the topology in $M$ (inherited from $I ^{\infty}$).

We denote by $C (I^{\infty} , I^{\infty} )$ the space of continuous maps of $I^{\infty}$ into itself, equipped with the compact - open topology.

Let us fix countable ordinals $\alpha$, $\beta$ and an $n$ which is a natural number or $\infty$.

Let $\Phi _{\alpha} : \bbN ^{ \bbN } \to \K _{\alpha} (M)$, $\Phi _{\beta} : \bbN ^{ \bbN } \to \K _{\beta} (M)$ be continuous surjections on the irrationals, considered in \cite{PP}, proof of Proposition 5.1.1. Recall that for every $\alpha < \omega _{1}$, the set $G_{\alpha}  = \{ (t,x): \, x \in \Phi _{\alpha} (t) \}$  is closed in $\bbN ^{ \bbN } \times M$, ind$G_{\alpha} = \alpha$ and
each $K \in \K _{\alpha} (M)$ is isometric with $\{t \} \times \Phi _{\alpha} (t)$, where $\Phi _{\alpha} (t) = K$.

Let
\smallskip\begin{enumerate}
\item[(1)] $ \G = \{ (s,t,f) \in \bbN ^{ \bbN } \times \bbN ^{ \bbN } \times C (I^{\infty} , I^{\infty} ):  f(\Phi _{\alpha}(s)) =  \Phi _{\beta}(t) \hbox{ and } {\rm dim}(f^{-1} (y) \cap \Phi _{\alpha}(s)) \leq n \hbox{ for each } y \in \Phi _{\beta}(t) \}$
\end{enumerate}
Using continuity of $\Phi _{\alpha}$ and $\Phi _{\beta}$, and the fact that at most $n$-dimensional compacta form a $G_{\delta}$-set in the hyperspace $\K (I^{\infty})$, cf. \cite{K}, $\S$45, IV, Theorem 4, we will show that
\smallskip\begin{enumerate}
\item[(2)] $\;\G$ is a $G_{\delta}$-set in the product $\bbN ^{ \bbN } \times \bbN ^{ \bbN } \times C (I^{\infty} , I^{\infty} )$.
\end{enumerate}
 Let $\E = \{ K \in \K (I ^{\infty}): \; {\rm dim}K > n \} $. Then $\E = \bigcup _{m} \E _{m}$, where $\E _{m}$ is closed in the hyperspace. Moreover, the set $\E _{m} ^{\star} = \{ L \in \K (I ^{\infty}): \; L \supset K \hbox{ for some } K \in \E _{m} \}$ is also closed in $\K (I ^{\infty})$, contains $\E _{m}$ and is contained in $\E$. Therefore, replacing $\E _{m}$ by $\E _{m} ^{\star}$
 we can assume that $\E _{m}$ contains all compacta containing some element of $\E _{m}$.

 Let us check that, for every $m \in \bbN$,
 \smallskip\begin{enumerate}
\item[(3)] $\; \F _{m} = \{  (s,f,y) \in  \bbN ^{ \bbN }\times C (I^{\infty} , I^{\infty} ) \times I^{\infty} : \; f^{-1} (y) \cap \Phi _{\alpha}(s) \in \E _{m} \} $ is closed in the product
$\bbN ^{ \bbN }\times C (I^{\infty} , I^{\infty} ) \times I^{\infty}$
\end{enumerate}
To that end, consider $ (s_{i}, f_{i}, y_{i}) \to (s,f,y)$,  $(s_{i}, f_{i}, y_{i}) \in \F _{m}$, and let $C_{i} = f_{i}^{-1} (y_{i}) \cap \Phi _{\alpha}(s_{i})  \in \E _{m}$.
Passing, if necessary, to a subsequence, we can assume that $C_{i} \to C$ in $\K (I^{\infty} )$.

We shall check that $C \subset f^{-1}(y) \cap \Phi _{\alpha}(s)$. Let us pick any $a \in C$, and let $a_{i} \to a$, $a_{i} \in C_{i}$. Since $f_{i}(a_{i}) = y_{i}$, $f_{i} \to f$, $y_{i} \to y$, we have
$f_{i} (a_{i}) \to f(a)$, hence $f(a)=y$. Also, since $a_{i} \in  \Phi _{\alpha}(s_{i})$ and $ \Phi _{\alpha}$ is continuous, $a \in  \Phi _{\alpha}(s)$. In effect, $a \in f^{-1}(y) \cap \Phi _{\alpha}(s)$, and $C \subset f^{-1}(y) \cap \Phi _{\alpha}(s)$. Since $C \in \E _{m}$ its superset $f^{-1}(y) \cap \Phi _{\alpha}(s)$ is also in $\E _{m}$, hence $(s,f,y) \in \F _{m}$, which gives us (3).

 Now, denoting by $\pi : \bbN ^{ \bbN }\times C (I^{\infty} , I^{\infty} ) \times I^{\infty}  \to \bbN ^{ \bbN }\times C (I^{\infty} , I^{\infty} )$ the projection parallel to the compact axis $I^{\infty}$, we conclude that the sets $\pi (\F _{m})$ are closed, and hence the sets
\smallskip\begin{enumerate}
\item[(4)] $\; \H _{m} =\left( \bbN ^{ \bbN }\times C (I^{\infty} , I^{\infty} )\right) \setminus \pi (\F _{m})$
\end{enumerate}
are open in $\bbN ^{ \bbN }\times C (I^{\infty} , I^{\infty} )$. Therefore, the set
\smallskip\begin{enumerate}
\item[(5)] $\; \bigcap _{m} \H _{m} = \{ (s,f) \in \bbN ^{ \bbN }\times C (I^{\infty} , I^{\infty} ): \; \hbox{for every } y \in I^{\infty}, \hbox{ dim}(f^{-1} (y) \cap  \Phi _{\alpha}(s)) \leq n \}$
\end{enumerate}
is a $G_{\delta}$-set in $\bbN ^{ \bbN }\times C (I^{\infty} , I^{\infty} )$.

Since $\L = \{ (s,t,f): \; f( \Phi _{\alpha}(s)) =  \Phi _{\beta}(t) \} \subset \bbN ^{ \bbN }\times \bbN ^{ \bbN }\times C (I^{\infty} , I^{\infty} ) $ is closed, both $\Phi _{\alpha}$, $\Phi _{\beta}$ being continuous, the set $\G = \L \cap \{ (s,t,f): \; (s,f) \in  \bigcap _{m} \H _{m} \}$ is a $G_{\delta}$-set in $\bbN ^{ \bbN } \times \bbN ^{ \bbN } \times C (I^{\infty} , I^{\infty} )$, which ends the proof of (2).

 By (2), there is a continuous surjection of $\bbN ^{ \bbN }$ onto $\G$,
\smallskip\begin{enumerate}
\item[(6)] $\; u \to  (s(u), t(u), f_{u} ) \in \G$, $\, u \in  \bbN ^{ \bbN }$,
\end{enumerate}
and we let
\smallskip\begin{enumerate}
\item[(7)] $\; X = \{ (u,z): \, z \in \Phi _{\alpha}(s(u)) \} \subset   \bbN ^{ \bbN }\times M$,
\end{enumerate}
\smallskip\begin{enumerate}
\item[(8)] $\; Y = \{ (t,z): \, z \in \Phi _{\beta}(t) \} \subset   \bbN ^{ \bbN }\times M$,
\end{enumerate}
\smallskip\begin{enumerate}
\item[(9)] $\; F : X \to Y$, $\, F(u,z) =  ( t(u), f_{u}(z) ) $.
\end{enumerate}

The spaces $X$, $Y$ are considered with the metric inherited from the product $ \bbN ^{ \bbN }\times M$, where $M$ is equipped with the universal metric and $ \bbN ^{ \bbN }$ has a standard complete metric.

Then (in notation from \cite{PP}), $Y = G_{\beta}$, hence ind$Y = \beta$, and since the function $(u,z) \to (s(u), z)$ maps $X$ into $G_{\alpha}$, using an observation in \cite{P1}, \S 4, Sublemma 3.2, we check that also ind$X = \alpha$.

The mapping $F$ in (9) is a continuous surjection. For any $(r,w)\in Y$, $F^{-1}(r,w) = \{ (u,z) \in \bbN ^{ \bbN }\times M: \; t(u)=r, \; z \in f_{u}^{-1}(w) \cap \Phi _{\alpha} (s(u))\}$, hence by (1) the projection of $F^{-1}(r,w)$ onto the first coordinate is a perfect map with at most $n$-dimensional fibers. Therefore, dim$F \leq n$.

Finally, let $f : K \to L$ be a continuous between compact spaces with ind$K \leq \alpha$, ind$L \leq \beta$ and dim$f \leq n$. Since $M$ was an isometrically universal space, we can assume that $K$ and $L$ are isometrically embedded in $M$. The function $f : K \to L$ can be extended to a continuous map $\tilde{f} : I^{\infty} \to I^{\infty}$. For $s,t \in \bbN ^{ \bbN }$ such that $K =  \Phi _{\alpha} (s)$, $L =  \Phi _{\beta} (t)$, we have $(s,t,\tilde{f}) \in \G$, cf. (1), and let $u \in \bbN ^{ \bbN }$  be such that $s=s(u)$, $t =t(u)$, $\tilde{f} = \tilde{f} _{u}$, cf. (6).

Then, for the isometric identifications $i : K \to \{ u \} \times K$ and $j : L \to \{t \} \times L$, we have $F \circ i = j \circ f$, cf. (9).

\bigskip

\section{Comment}

In \cite{O}, Section 4, A.I.Oblakova proved that there exists a Cantor set such that any finite metric space whose diameter does not exceed $1$ and the number of points does not exceed $n$ can be isometrically embedded into it.
 Using a continuous parametrization of some collections of finite metric spaces, one can   refine slightly this result, to the following effect:

\bigskip

{\bf Remark.} {\it For each natural number $n$ there are Cantor sets  $X$ and $Y$ and a continuous surjection $F : X \to Y$ which is isometrically containing for the class $\D _{n}$ of non-expanding surjections between at most $n$-element metric spaces of diameter $\leq 1$.}

\bigskip

\begin{proof} Let us fix a natural number $n$. First we will prove that
\smallskip\begin{enumerate}
\item[(10)] there are zero-dimensional compact metric spaces $E,H$ and a continuous surjection $f : E \to H$ which is isometrically containing for the class $\D_{n}$.
\end{enumerate}
Indeed, by a result of Gromov \cite{G}, sec.6, there is a compact metric space $Z$ which contains isometrically each metric space of diameter $\leq 1$ containing at most $n$ elements (cf. also \cite{O} and \cite{I2}, 9.3, for different proofs).

Let $Z \times Z$ be the metric product and $\pi _{i} : Z \times Z \to Z$, $i=1,2$, projections onto the first and the second coordinate, respectively.

Let $\F$ be the collection of at most $n$-element subsets $T$ of $Z \times Z$ which are graphs of non-expanding surjections from $\pi _{1} (T)$ onto $\pi _{2} (T)$, cf. \cite{PP}, (4).

Then $\F$ is compact in the hyperspace of $Z \times Z$ and let $u : 2^{\bbN} \to \F$ be a continuous surjection of the Cantor set $2^{\bbN}$ onto $\F$.

We let
\smallskip\begin{enumerate}
\item[(11)] $\;E= \{ (t,x): \; x \in \pi _{1} (u(t)) \} \subset  2^{\bbN} \times Z$,
\end{enumerate}
\smallskip\begin{enumerate}
\item[(12)] $H= \{ (t,y): \; y \in \pi _{2} (u(t)) \} \subset  2^{\bbN} \times Z$,
\end{enumerate}
\smallskip\begin{enumerate}
\item[(13)] $f : E \to H$, $\; f(t,x) = (t,y)$, where $(x,y) \in u(t)$.
\end{enumerate}
\noindent Notice that continuity of $f$ follows from the fact that its graph $\{ \left( (t,x), (t,y) \right) : \; (x,y) \in u(t) \}$ is a closed subset of the compact product
$(2^{\bbN} \times Z) \times (2^{\bbN} \times Z)$.

To see that $f$ is isometrically containing for the class $\D _{n}$, suppose that $g : K \to L$ is a non-expanding surjection between at most $n$-element metric spaces of diameter $\leq 1$.
One can assume that $K$ and $L$ are isometrically embedded in $Z$.  Since $T = \{ (x,g(x)): x \in K \}$ belongs to $\F$, $T = u(t)$ for some $t \in 2^{\bbN}$. Then, for the isometric identifications $i : K \to \{ t \} \times K$ and $j : L \to \{t \} \times L$, we have $f \circ i = j \circ g$, cf. (13).

To end the proof of Remark, it suffices to embed $E$ and $H$ into Cantor sets $X$ and $Y$ (by a theorem of Hausdorff, one can extend the metrics on $E$ and $H$ to metrics on $X$ and $Y$, respectively), and to define $F : X \to Y$ as $F = f \circ r$, where $r : X \to E$ is a retraction. Finally, to make sure that $F$ is a surjection, we can always add to $X$ a disjoint copy of $Y$ and let $F$ be the identity on this copy.
\end{proof}

\end{document}